\documentclass[reqno,draft]{amsart}

\input cyracc.def

\newcommand{\labell}[1]{\label{#1}}

\DeclareMathOperator{\Ad}{Ad}

\newcommand{\myvec}[1]{{\text{\bf #1}}}

\newtheorem{thm}{Theorem}
\newtheorem*{thm*}{Theorem}
\theoremstyle{remark}
\newtheorem*{rmk*}{Remark}
\newtheorem*{ack*}{Acknowledgements}

\newcounter{cond}

\begin{document}

\title[Equivariant Liapunov Stability Test and
Energy--Momentum--Casimir]{An Equivariant Liapunov Stability Test and
the Energy--Momentum--Casimir Method} \author{Eric T. Matsui} 
\address{University of California, Santa Cruz}
\email{ematsui@math.ucsc.edu}
\date{February 11, 2002}
\begin{abstract}
  We present an equivariant Liapunov stability criterion for dynamical
  systems with symmetry.  This result yields a simple proof of the
  energy--momentum--Casimir stability analysis of relative equilibria 
  of equivariant Hamiltonian systems.
\end{abstract}
\maketitle

\section{Introduction}
We present an equivariant Liapunov stability criterion for dynamical
systems with symmetry.  This result yields a simple proof of the
energy--momentum--Casimir (EMC for short) stability analysis of
relative equilibria of equivariant Hamiltonian systems.

The equivariant Liapunov stability test can be viewed as the natural
adaptation of the classic Liapunov stability theorem to relative
equilibria of equivariant dynamical systems.  The proof echoes the
traditional proof of Liapunov stability with all objects replaced by
their equivariant/invariant counterparts.

The Liapunov stability theorem has numerous applications to
Hamiltonian systems, where the Hamiltonian and other conserved
quantities provide the ingredients of the Liapunov functions.

The EMC method is widely used in applications to obtain sufficient
conditions for group-stability of relative equilibria.  The two key
components of the EMC method are the construction of a suitable
augmentation of the Hamiltonian and the application of a second
derivative test (or, in the infinite-dimensional case, a convexity
estimate) to establish stability.

The EMC method has been developed in a series of papers.  The earliest
version yielding rigorous group stability results used convexity
arguments for systems on vector spaces; see \cite{cit:hmrw85} and
references therein.  For free actions on manifolds, the result has
been obtained by Patrick, \cite{cit:patrick92}.  Non-free actions
(with \( G_{\mu} \) abelian) have been treated by Lewis,
\cite{cit:lew92}.  The most general result to date (proper actions, \(
G_{\mu} \) compact) is due to Ortega, \cite{cit:ortega98}.

In Section~\ref{sec:emc} we show that the application of the
equivariant Liapunov stability test to this augmented Hamiltonian
immediately yields Ortega's result.  This shows that the only aspect
of the Hamiltonian structure utilized in the stability analysis is the
existence of an invariant conserved quantity.  In particular, the 
analysis does not make use of the symplectic or Poisson structure.  

In the last section, we present an extension of our results to
infinite-dimensional systems.  We adapt the key ideas behind Arnol'd's
method of convexity estimates (\cite{cit:arn69,cit:hmrw85}) to the
group-equivariant setting.

These results comprise a part of the author's Ph.D. thesis 
(\cite{cit:matsui01a}).  An application of the EMC method to the 
stability and bifurcation analysis of the sleeping pseudo-Lagrange top 
also appears in \cite{cit:matsui01a}.  This application is further 
developed in \cite{cit:matsui01b} including an interesting example of 
a symmetry-breaking bifurcation to a family of steadily precessing 
motions where no ``exchange of stability'' occurs; both symmetric and 
asymmetric branches are stable in the neighborhood of the bifurcation 
point.
\section{Equivariant Liapunov Stability}
\labell{sec:EL}
Let \( X \) be a smooth vector field on a smooth finite-dimensional
manifold \( M \) that is equivariant with respect to a smooth and
proper action of some Lie group \( G \).  Let \( \mathfrak{g} \)
denote the Lie algebra of \( G \), and let \( F_{t} \), \(
t\in\mathbb{R} \), denote the flow of \( X \).

An element \( m\in M \) is called a \emph{relative equilibrium} if the
curve \( (\exp t\xi)\cdot m \) is a trajectory of \( X \) for some \(
\xi \) in \( \mathfrak{g} \).

Let \( G' \) be a subgroup of \( G \).  A relative equilibrium \( m \)
of \( X \) is said to be \emph{\( G' \)-stable}, or \emph{stable
modulo \( G' \)}, if for every \( G' \)-invariant open neighborhood \(
U \) of \( G'\cdot m \) there exists an open neighborhood \( V \) of \(
m \) in \( U \) such that \( F_{t}(V)\subseteq U \) for all \( t\geq 0
\).

\begin{thm}[The Equivariant Liapunov Stability Test]
  \labell{thm:els}
  Let \( m\in M \) be a relative equilibrium.  If there exists a \( G
  \)-invariant function \( f \) such that
  \begin{list}{\rm LS\arabic{cond}.}{\usecounter{cond}}
  \item\labell{cond:lscp} \( m \) is a critical point of \( f \) (say,
    \( f(m)=0 \)),
  
  \item\labell{cond:lssemipd} the bilinear form \( D^{2}f(m) \) is
    positive semi-definite with kernel \( \mathfrak{g}\cdot m \), and
  
  \item\labell{cond:cnstrnt} for every \( \epsilon > 0 \) there exists
    an open neighborhood \( V_{\epsilon} \) of \( m \) such that
    \begin{equation*}
      f(F_{t}(V_{\epsilon})) \subset[0,\epsilon),\quad\forall\ t>0,
    \end{equation*}
  \end{list}
  then \( m \) is \( G \)-stable.
\end{thm}

\begin{proof}
  Without loss of generality, we may assume that \( f(m)=0 \).
  
  Let \( U \) be a \( G \)-invariant open neighborhood of the orbit \(
  G\cdot m \).  Our task is to construct an open neighborhood \( V \)
  of \( m \) such that \( F_{t}(V)\subseteq U \) for all \( t>0 \).
  
  By the Slice Theorem (see, e.g., \cite{cit:ggk}), there exists a \(
  G \)-invariant open neighborhood \( \tilde{M} \) of \( G\cdot m \)
  and a \( G \)-equivariant projection map \(
  \pi\colon\tilde{M}\rightarrow G\cdot m \) that makes \( \tilde{M} \)
  a disc bundle.  Let \( W \) denote the fiber \( \pi^{-1}(m) \).
  
  Properties LS\ref{cond:lscp}, LS\ref{cond:lssemipd}, and the Morse
  Lemma guarantee the existence of positive number \( \epsilon \) such
  that the connected component \( W_{\epsilon} \) of \(
  \left(\left.f\right|_{W}\right)^{-1}[0,\epsilon) \) is an open disk
  embedded in \( W \) and contained in \( U \).  Let \( G^{0} \)
  denote the connected component of \( G \) containing the identity. 
  Since \( U \) is \( G \)-invariant, \( G^{0}\cdot W_{\epsilon} \) is
  contained in \( U \).  Using the associated bundle structure of \(
  \tilde{M} \) (see \cite{cit:ggk}), one can show that \( G^{0}\cdot
  W_{\epsilon} \) is open and is the connected component of \(
  f^{-1}[0,\epsilon) \) containing \( m \).
  
%
  Let \( V=(G^{0}\cdot W_{\epsilon})\cap V_{\epsilon} \), and let \( v
  \) be an element of \( V \).  Since \( m \) is an element of \( V
  \), LS\ref{cond:cnstrnt} implies that \( F_{t}(v) \) is in \(
  [0,\epsilon) \) for all \( t>0 \), so that \( F_{t}(v) \) is in the
  connected component of  \( f^{-1}[0,\epsilon) \)
 containing \( v \), and, hence, \( m \). 
  This implies that \( F_{t}(v) \) is in \(
  G^{0}\cdot W_{\epsilon} \) for all \( t>0 \).  Since \( G^{0}\cdot 
  W_{\epsilon}\subseteq U \), this proves that \( m \) is \( G 
  \)-stable.
\end{proof}
\section{Application to Hamiltonian Systems with Symmetry}
\labell{sec:emc}
Let \( (M,\{\cdot,\,\cdot\}) \) be a finite-dimensional Poisson
manifold.  Let \( G \) act on \( M \) properly and canonically.  Let
us assume that there exists an equivariant momentum map \(
\myvec{J}\colon M\rightarrow\mathfrak{g}^{*} \) associated to the \( G
\)-action.  Let \( H \) be a \( G \)-invariant Hamiltonian and \(
X_{H} \) its Hamiltonian vector field.  Let \( z_{e}\in M \) be a
relative equilibrium of \( X_{H} \) satisfying \(
X_{H}(z_{e})=\xi_{M}(z_{e}) \), where \( \xi_{M} \) denotes the
infinitesimal action on \( M \) by some Lie algebra element \(
\xi\in\mathfrak{g} \).  Set \( \mu=\myvec{J}(z_{e}) \).

Let us assume that there exist \( G_{\mu} \)-invariant norms \( 
\lVert\cdot\rVert_{\mathfrak{g}} \) and \( 
\lVert\cdot\rVert_{\mathfrak{g}^{*}} \) on \( \mathfrak{g} \) and \( 
\mathfrak{g}^{*} \), respectively, that satisfy
\begin{equation}
  \langle\alpha,\zeta\rangle
  \leq\lVert\alpha\rVert_{\mathfrak{g}^{*}}\,
  \lVert\zeta\rVert_{\mathfrak{g}}
  \labell{eq:csineq}
\end{equation}
for all \( \alpha\in\mathfrak{g}^{*} \), \( \zeta\in\mathfrak{g} \). 
For example, this is the case if \( G_{\mu} \) is compact.  Let \( V
\) be an inner product space, and let \( \myvec{C}:M\rightarrow V \)
be a map that is invariant with respect to the \( G_{\mu} \)-action
and the flow.  An energy--momentum--Casimir function is defined by
\begin{displaymath}
  \mathcal{H}(z):=H(z)-J_{\xi}(z)+C_{\lambda}(z),\quad z\in M,
\end{displaymath}
  where \( J_{\xi}(z):=\langle\myvec{J}(z),\xi\rangle \) and \(
  C_{\lambda}(z):=\langle\lambda,\myvec{C}(z)\rangle \) for some \(
  \lambda\in V^{*} \).
\begin{thm}[Energy--Momentum--Casimir Method,
  \cite{cit:lew92, cit:ortega98, cit:patrick92}]
  \labell{thm:emc}
  Suppose that there
  exists a \( \lambda \) in \( V^{*} \) such that the following
  conditions hold:
  \begin{list}{\rm EM\arabic{cond}.}{\usecounter{cond}}
    \item\labell{cond:emcp} \( z_{e} \) is a critical point of \( 
    \mathcal{H} \).
  
    \item\labell{cond:emsemipd} The restriction of the bilinear form
    \(
    D^{2}\mathcal{H}(z_{e})
    \) to
    \begin{equation}
      \labell{eq:k}
      K=\ker D\myvec{J}(z_{e})\cap\ker D\myvec{C}(z_{e})
    \end{equation}
    is positive or negative semi-definite with
    kernel \( \mathfrak{g}_{\mu}\cdot z_{e} \).
  
    \item\labell{cond:gzingxi} \( G_{z_{e}}\subseteq G_{\xi} 
    \).
  \end{list}
  Then \( z_{e} \) is \( G_{\mu} \)-stable.
\end{thm}
\begin{rmk*}
  The condition EM\ref{cond:gzingxi} is trivially satisfied if the \(
  G \)-action is free, since \( G_{z_{e}} \) would be trivial.  If \(
  G_{z_{e}} \) is nontrivial, then \( \xi \) satisfying \(
  X_{H}(z_{e})=\xi_{M}(z_{e}) \) is determined up to \(
  \mathfrak{g}_{z_{e}} \).  Then \( \xi \) can always be chosen so
  that EM\ref{cond:gzingxi} is satisfied.  The argument is as follows.
  
  Since \( G_{z_{e}} \) is compact, there exists a \( G_{z_{e}}
  \)-invariant complement \( \mathfrak{g}_{z_{e}}^{\perp} \) to \(
  \mathfrak{g}_{z_{e}} \) in \( \mathfrak{g} \).  Let \( \xi^{\perp}
  \) be the \( \mathfrak{g}_{z_{e}}^{\perp} \)-component of \( \xi \). 
  By the equivariance of \( X_{H} \), we have \( g\cdot
  X_{H}(z_{e})=X_{H}(g\cdot z_{e})=X_{H}(z_{e}) \) for \( g \) in \(
  G_{z_{e}} \).  Since \(
  X_{H}(z_{e})=\xi_{M}(z_{e})={\xi^{\perp}}_{M}(z_{e}) \), this
  implies that \( g\cdot{\xi^{\perp}}_{M}(z_{e})=
  (\Ad_{g}\xi^{\perp})_{M}(z_{e})={\xi^{\perp}}_{M}(z_{e}) \).  Hence,
  \( \Ad_{g}\xi^{\perp}-\xi^{\perp} \) is in \( \mathfrak{g}_{z_{e}}
  \).  On the other hand, \( \mathfrak{g}_{z_{e}}^{\perp} \) is \(
  G_{z_{e}} \)-invariant, so that \( \Ad_{g}\xi^{\perp}-\xi^{\perp} \)
  is in \( \mathfrak{g}_{z_{e}}^{\perp} \).  Therefore, \(
  \Ad_{g}\xi^{\perp}=\xi^{\perp} \), and \( G_{z_{e}}\subseteq
  G_{\xi^{\perp}} \).  The algebra element \( \xi \) satisfying \(
  X_{H}(z_{e})=\xi_{M}(z_{e}) \) and EM\ref{cond:gzingxi} is not
  necessarily unique.  This is important, because the difference in
  the choice of \( \xi \) influences the validity of
  EM\ref{cond:emsemipd}.  The relative equilibrium is stable as long
  as EM1--3 are satisfied for at least one choice of \( \xi \), and
  having more choices increases the chance of establishing stability.
  
  A case when EM\ref{cond:gzingxi} is valid for all \( \xi \)
  satisfying \( X_{H}(z_{e})=\xi_{M}(z_{e}) \) is when \( G_{\mu} \)
  is abelian, which is generically the case if \( G_{\mu} \) is
  connected (see \cite{cit:dv69, cit:mr94}).  In this case, since \(
  \xi \) is an element of \( G_{\mu} \), it follows that \(
  G_{\mu}\subseteq G_{\xi} \).  Since \( G_{z_{e}} \) is contained in
  \( G_{\mu} \), EM\ref{cond:gzingxi} is satisfied.
\end{rmk*}
\begin{proof}
  We first prove the case when \(
  D^{2}\mathcal{H}(z_{e}) \) is positive
  semi-definite on \( K \).
  
  The condition EM\ref{cond:gzingxi} is necessary for the following
  construction.  By the Slice Theorem, there exists a \( G_{\mu}
  \)-invariant open neighborhood \( \tilde{M} \) of \( z_{e} \) and a
  \( G_{\mu} \)-equivariant projection map \(
  \pi\colon\tilde{M}\rightarrow G_{\mu}\cdot z_{e} \) that makes \(
  \tilde{M} \) a vector bundle.  Let \( \phi\colon
  G_{\mu}/G_{z_{e}}\rightarrow G_{\mu}\cdot z_{e} \) denote the
  natural isomorphism.  Define \(
  \tilde{\pi}\colon\tilde{M}\rightarrow G_{\mu}/G_{z_{e}} \) as \(
  \tilde{\pi}=\phi^{-1}\circ\pi \).  The condition
  EM\ref{cond:gzingxi} guarantees that
  \begin{displaymath}
    \Ad_{\tilde{\pi}(z)}\xi,\quad z\in M,
  \end{displaymath}
  is well defined.  This is called \emph{Patrick's velocity map}
  (\cite{cit:ortega98}).  We incorporate it in the definition of the
  \( G_{\mu} \)-invariant function used in the proof of the
  equivariant Liapunov stability test.  Let
  \begin{displaymath}
    f(z)=f_{1}(z)+\sigma f_{2}(z),\quad z\in M,
  \end{displaymath}
  where \( \sigma \) is a real number sufficiently large so that \( f
  \) satisfies LS\ref{cond:lssemipd} (see below), and \( f_{1} \) and
  \( f_{2} \) are defined by
  \begin{align*}
    f_{1}(z)&=H(z)-\langle\myvec{J}(z),\cdot\Ad_{\tilde{\pi}(z)}\xi
    \rangle
    +C_{\lambda}(z)\\
    f_{2}(z)&=\lVert\myvec{J}(z)-\mu\rVert_{\mathfrak{g}^{*}}^{2}+
    \lVert \myvec{C}(z)-\myvec{C}(z_{e})\rVert_{V}^{2}.
  \end{align*}
  The \( G_{\mu} \)-equivariance of \( \myvec{J} \) and \( \tilde{\pi}
  \) guarantees that \( f_{1} \) is \( G_{\mu} \)-invariant.  The \(
  G_{\mu} \)-invariance of \( f_{2} \) follows from the equivariance
  of \( \myvec{J} \) and \( G_{\mu} \)-invariance of \(
  \lVert\cdot\rVert_{\mathfrak{g}^{*}} \) and \( \myvec{C} \).  Hence,
  \( f \) is \( G_{\mu} \)-invariant.
  
  We now prove that \( f \) satisfies conditions LS\ref{cond:lscp},
  LS\ref{cond:lssemipd}, and LS\ref{cond:cnstrnt}, for sufficiently 
  large values of \( \sigma \).
  
  \emph{Condition LS\ref{cond:lscp}}.\quad That \( z_{e} \) is a
  critical point of \( f_{2} \) is obvious.  Hence, it suffices to
  prove that \( z_{e} \) is a critical point of \( f_{1} \).
  
  Let \( S \) denote the slice \( \pi^{-1}(z_{e}) \).  Then
  \begin{equation}
    T_{z_{e}}M=\left(\mathfrak{g}_{\mu}\cdot z_{e}\right)
    \oplus\left(T_{z_{e}}S\right).  \labell{eq:slicedecomp}
  \end{equation}
  Since \( f_{1} \) is \( G_{\mu} \)-invariant, this implies that we
  only need to show that \( z_{e} \) is a critical point of \(
  f_{1}|_{S} \).  Since \(
  \tilde{\pi}(S)=\phi^{-1}(z_{e})=G_{z_{e}}\subseteq G_{\xi} \) (where
  \( \phi:G_{\mu}/G_{z_{e}}\rightarrow G_{\mu}\cdot z_{e} \) is the
  natural isomorphism) we have \( \Ad_{\tilde{\pi}(z)}\xi=\xi \) for
  \( z \) in \( S \), so that
  \begin{equation}
    f_{1}(z)=\mathcal{H}(z),\quad \text{for\
    }z\in S. \labell{eq:Hxiandf1onS}
  \end{equation}
  Therefore, by EM\ref{cond:emcp}, \( z_{e} \) is a critical point of
  \( f_{1}|_{S} \).  This proves that \( z_{e} \) is a critical point
  of \( f=f_{1}+\sigma f_{2} \).

  \emph{Condition LS\ref{cond:lssemipd}}.\quad From the decomposition
  \eqref{eq:slicedecomp}, we see that the condition
  LS\ref{cond:lssemipd} is satisfied if and only if \( D^{2}f(z_{e})
  \) is positive-definite on \( T_{z_{e}}S \).  Let \(
  D_{S}^{2}f_{1}(z_{e}) \) and \( D_{S}^{2}f_{2}(z_{e}) \) denote the
  restrictions of \( D^{2}f_{1}(z_{e}) \) and \( D^{2}f_{2}(z_{e}) \)
  to \( T_{z_{e}}S \), respectively.  We claim that \(
  D_{S}^{2}f_{2}(z_{e}) \) is positive semi-definite with kernel \(
  K\cap T_{z_{e}}S \), and that \( D_{S}^{2}f_{1}(z_{e}) \) is
  positive-definite on \( K\cap T_{z_{e}}S \).  Then \(
  D_{S}^{2}f_{1}(z_{e})+\sigma D_{S}^{2}f_{2}(z_{e}) \) is
  positive-definite for sufficiently large values of \( \sigma \), so
  that \( f \) satisfies LS\ref{cond:lssemipd}.
  
  The claim is proved as follows.
  
  Since \( D^{2}f_{2}(z_{e}) \) is positive semi-definite with kernel
  \( K \), we conclude that the quadratic form \(
  D_{S}^{2}f_{2}(z_{e}) \) is positive semi-definite with kernel \(
  K\cap T_{z_{e}}S \).
  
  Since \( K \) contains \( \mathfrak{g}_{\mu}\cdot z_{e} \), the
  decomposition \eqref{eq:slicedecomp} implies that
  \begin{displaymath}
    K=\left(\mathfrak{g}_{\mu}\cdot z_{e}\right)\oplus\left(K\cap
    T_{z_{e}}S\right).
  \end{displaymath}
  Therefore, by EM\ref{cond:emsemipd}, \(
  D^{2}\mathcal{H}(z_{e}) \) is positive-definite on \(
  K\cap T_{z_{e}}S \).  Combining this with \eqref{eq:Hxiandf1onS}, we
  conclude that \( D_{S}^{2}f_{1}(z_{e}) \) is positive-definite on \(
  K\cap T_{z_{e}}S \).

  \emph{Condition LS\ref{cond:cnstrnt}.}\quad Let \( F_{t} \) denote
  the flow of \( X_{H} \).  Since \(
  \tilde{\pi}(z_{e})=G_{z_{e}}\subseteq G_{\xi} \), we have \(
  \langle\myvec{J}(z_{e}),\Ad_{\tilde{\pi}(z_{e})}\xi\rangle
  =\mu\cdot\xi \), and
  \begin{align*}
    f(F_{t}(z))-f(z_{e})&=H(F_{t}(z))-H(z_{e})-\left(\langle
    \myvec{J}(F_{t}(z)),
    \Ad_{\tilde{\pi}(F_{t}(z))}\xi\rangle-\langle\mu,\xi\rangle\right)\\
    &\quad +C_{\lambda}(F_{t}(z))-C_{\lambda}(z_{e})+\sigma 
    f_{2}(F_{t}(z))\\
    &=H(z)-H(z_{e})-\left(\langle\myvec{J}(z),\Ad_{\tilde{\pi}(F_{t}(z))}
    \xi\rangle-\langle\mu,\xi\rangle\right)\\
    &\quad +C_{\lambda}(z)-C_{\lambda}(z_{e})+\sigma f_{2}(z).
  \end{align*}
  Since \( \tilde{\pi}(z) \) is in \( G_{\mu}/G_{z_{e}} \) for any \(
  z\in \tilde{M} \), we have \( \mu=\Ad_{\tilde{\pi}(F_{t}(z))}^{*}\mu
  \), and
  \begin{align*}
    f(F_{t}(z))-f(z_{e})&=H(z)-H(z_{e})-\langle\myvec{J}(z)
    -\mu,\Ad_{\tilde{\pi}(F_{t}(z))}\xi\rangle\\
    &\quad +C_{\lambda}(z)-C_{\lambda}(z_{e})+\sigma f_{2}(z)\\
    &\leq\left\lvert\,H(z)-H(z_{e})\,\right\rvert+\lVert\,
    \myvec{J}(z)-\myvec{J}(z_{e})\,\rVert_{\mathfrak{g}^{*}}\,
    \lVert\,\xi\,\rVert_{\mathfrak{g}}\\
    &\quad +\lvert C_{\lambda}(z)-C_{\lambda}(z_{e})\rvert+\sigma
    \lvert f_{2}(z)\rvert.
  \end{align*}
  The inequality follows from \eqref{eq:csineq} and \( G_{\mu}
  \)-invariance of \( \lVert\,\cdot\,\rVert_{\mathfrak{g}} \).  Hence,
  LS\ref{cond:cnstrnt} is satisfied by the continuity of \( H \), \(
  \myvec{J} \), and \( \myvec{C} \).

  If \( D^{2}\mathcal{H}(z_{e}) \) is negative semi-definite on \( K
  \) with kernel \( \mathfrak{g}_{\mu}\cdot z_{e} \), then we redefine
  \( f \) by
  \begin{displaymath}
    f(z)=-f_{1}(z)+\sigma f_{2}(z),\quad z\in M.
  \end{displaymath}
  The proof proceeds as in the positive semi-definite case.
\end{proof}


\section{Extension to Infinite Dimensions} It is well-known that the
positive-definiteness of the Hessian of the energy at an equilibrium
is a sufficient condition for its stability, provided that the system
is finite-dimensional.  For infinite-dimensional systems, additional
assumptions are necessary to establish stability.  The convexity
estimate method due to Arnol'd has been widely used to circumvent this
shortcoming in order to establish nonlinear stability.  (See,
\textit{e.g.}, \cite{cit:hmrw85}, and references therein.)  In the
remainder of this article, we present a possible adaptation of
Arnol'd's approach to the stability of relative equilibria with
isotropy in infinite dimensions.  We make use of the Slice Theorem for
the action of compact groups on Banach manifolds, which allows us to
easily extend the results of the previous two sections to this
setting.  (See, \textit{e.g.}, \cite{cit:ggk} for the precise
statement of the Slice Theorem on Banach manifolds.)

\subsection{Equivariant Liapunov Stability}
Let \( M \) be a Banach manifold and let \( G \) be a compact Lie
group acting on \( M \).  Let \( X \) be a \( G \)-equivariant vector
field on \( M \).  Note that we do not assume that \( X \) is
Lipschitz or even continuous.  Relative equilibria of \( X \) can
still be defined in the usual way.  (However, without suitable
smoothness assumptions on \( X \) we cannot guarantee existence and
uniqueness of an integral curve for any initial conditions.)

\begin{thm}
  \labell{thm:bEL} Let \( m \) be a relative equilibrium.  Fix the
  Slice Theorem decomposition \( G\times_{G_m}S \) of a neighborhood
  of \( G\cdot m \), where \( S \) is an open ball in a Banach space.
  Without loss of generality, we may assume that the norm on \( S \)
  is \( G_m \)-invariant.  Let \( P:M\rightarrow\mathbb{R} \) be a
  continuous \( G \)-invariant function whose restriction to \( S \)
  is a function of the norm which has a strict minimum (say, equal to
  zero) at the origin.  If there exists a \( G \)-invariant function
  \( f \), with \( f(m)=0 \), such that
  \begin{list}{\rm LS\arabic{cond}$^\prime$.}{\usecounter{cond}}
  \item \( m \) is a critical point of \( f \),
  \item \labell{itm:Mbound} \( P(u)\leq f(u) \) for all \( u \) near
    \( m \), and
  \item for every \( \epsilon>0 \) there exists a neighborhood \( V \)
    of \( m \) such that \( f(u(t))<\epsilon \) for any integral curve
    \( u(t) \) with initial condition \( u(0) \) in \( V \).
  \end{list}
  Then for all \( \epsilon>0 \) there exists a neighborhood \( U \) of
  the orbit such that \( P(u(t))<\epsilon \) for any integral curve \(
  u(t) \) with \( u(0) \) in \( U \) as long as \( u(t) \) is defined.
\end{thm}
The proof of this theorem is essentially identical to the proof of the
Theorem~\ref{thm:els}.  The role of condition
LS\ref{itm:Mbound}$^\prime$ is to ensure that \( f^{-1}[0,\epsilon) \)
contains a small neighborhood of the orbit \( G\cdot m \). (In the
finite-dimensional version, this follows from LS\ref{cond:lssemipd}
and the Morse Lemma.)
\begin{rmk*}\ 
  \begin{enumerate}
  \item By the Slice Theorem, \( P \) only needs to be defined on the
    slice through \( m \).
  \item If \( P \) is the square of the original Banach norm on the
    slice, we recover the \( G \)-stability as in \S\ref{sec:EL}.
    Sometimes, \( P \) is not even the square of a norm; see
    \cite{cit:lew93a,cit:hes80}.
  \item With suitable modifications, the theorem can be extended to
    the case where \( P \) is a function of a norm which is weaker
    than the original norm on \( S \).  This generalization is useful
    when the original norm is chosen to be as close as possible to
    guaranteeing the existence and uniqueness of solutions and the
    ``norm square'' \( P \) is defined by the Hessian of \( f \) (cf.
    \cite{cit:bru02,cit:hmrw85}).
  \end{enumerate}
\end{rmk*}

\subsection{The EMC method}
Let \( M \) be a smooth Banach Poisson manifold and let \( G \) be a
compact Lie group that acts on it canonically.  In the following theorem,
we let \( S \) be the \( G_\mu \)-slice through \( z_e \).  By the
Slice Theorem, we can think of \( S \) as a ball, centered at the
origin, in a Banach space.  Let
\( \mathcal{H} \) and \( f_2 \) be defined as before.  Set \(
\mathcal{H}_\sigma =\mathcal{H}+\sigma f_2 \).
\begin{thm}
  \labell{thm:bEMC}
  Suppose that \( H \), \( \mathbf{J} \), and \( \mathbf{C} \) are continuous and there exists \( \lambda \in V^{*} \) and \(
  \sigma \in \mathbb{R} \) such that the following conditions
  hold:
  \begin{list}{\rm EM\arabic{cond}$^\prime$.}{\usecounter{cond}}
  \item \( z_{e}\in M \) is a critical point of \( 
    \mathcal{H} \).
  \item \labell{itm:convex} There exists a real number \( C \) such that
    \begin{displaymath}
      \lVert\Delta z\rVert^2\leq D^2\mathcal{H}_\sigma(z)(\Delta z,\Delta z)
    \end{displaymath}
    for all \( z \) in \( S \) near the origin and for all \( \Delta
    z\) in \( S \).
  \item \( G_{z_{e}}\subseteq G_{\xi} \).
  \end{list}
  Then \( z_e \) is Liapunov stable in the sense of
  Theorem~\ref{thm:bEL}.
\end{thm}
This theorem is proved by using Theorem~\ref{thm:bEL} in the same way
as Theorem~\ref{thm:els}.  Note that \( \mathcal{H}_\sigma \) and \(
f=f_1+\sigma f_2 \) (as defined in the proof of the finite-dimensional
EMC theorem) are identical on \( S \).  The difference is that \( f \)
is globally \( G_\mu \)-invariant, while \( \mathcal{H}_\sigma \) is
not.  Hence, we apply Theorem~\ref{thm:bEL} to \( f \), not \(
\mathcal{H}_\sigma \).  The condition LS\ref{itm:Mbound}$^\prime$
follows from EM\ref{itm:Mbound}$^\prime$ by the Mean Value Theorem.

\begin{ack*}
  The author would like to thank Viktor Ginzburg and Debra Lewis for
  their invaluable support during the preparation of this article.
\end{ack*}

\bibliographystyle{amsalpha}
\bibliography{bib}
\end{document}